\documentclass{mfatshort}
\usepackage{textcase}
\usepackage[english]{babel}
\usepackage[T2A]{fontenc}
\usepackage{amsmath}
\usepackage{amssymb}

\usepackage{amscd}
\usepackage{amsthm}
\usepackage{amstext}
\usepackage{amsfonts}
\usepackage{amsopn}
\usepackage{euscript}
\newcommand{\limind}{\mathop{{\rm lim\,ind}}\limits}
\newcommand{\XX}{{\mathfrak X}}
\newcommand{\N}{{\mathbb N}}

\newcommand{\R}{{\mathbb R}}
\newcommand{\CC}{{\mathbb C}}
\newcommand{\D}{{\mathcal D}}
\newcommand{\E}{{\mathfrak E}}
\newcommand{\I}{{\mathbb I}}
\newcommand{\Imm}{{\bf\rm Im}\,}
\newcommand{\K}{{\mathcal K}}
\newcommand{\QQ}{{\mathfrak Q}}
\newcommand{\cc}{\widetilde{c}}

\newcommand{\xx}{\widetilde{x}}

\newtheorem{theorem}{Theorem}[section]

\newtheorem{lemma}{Lemma}[section]
\newtheorem{nasl}{Corollary}[section]
\theoremstyle{remark}
\date{\empty}
\newtheorem{rmk}{Remark}[section]
\newtheorem{example}{Example}[section]

\begin{document}

\author{Ya. Grushka and S. Torba}

\title[Direct Theorems in the Theory of Approximation ...]{Direct Theorems in the Theory of
Approximation of Banach Space Vectors by Exponential Type Entire Vectors}

\email{sergiy.torba@gmail.com, grushka@imath.kiev.ua}
\address{Institute of Mathematics of the National Academy of
Sciences of Ukraine, Tereshchenkovskaya 3, 01601 Kiev (Ukraine)}

\thanks{This work was partially supported by the Ukrainian State
Foundation for Fundamental Research (project N14.1/003).}

\date{04/04/2007}

\begin{abstract}
For an arbitrary operator $A$ on a Banach space $\XX$ which is the
generator of $C_0$--group with certain growth condition at infinity,
the direct theorems on connection between the smoothness
degree of a vector $x\in\XX$ with respect to the operator $A$, the
rate of convergence to zero of the best approximation of $x$ by
exponential type entire vectors for the operator $A$, and the
$k$-module of continuity are established. The results allow to
obtain Jackson-type inequalities in a number of classic spaces of periodic
functions and weighted $L_p$ spaces.
\end{abstract}

\maketitle

\section{Introduction}

The direct and inverse theorems establishing a relationship
between the smoothness degree of a function with respect to the
differentiation operator and the rate of convergence to zero of
its best approximation by trigonometric polynomials are well known
in the theory of approximation of periodic functions. Jackson's
inequality is one among such results.

N. P. Kuptsov proposed a generalized notion of the module of
continuity, expanded onto $C_0$-groups in a Banach space
\cite{Kyptsov}. Using this notion, N. P.
Kuptsov \cite{Kyptsov} and A. P. Terekhin \cite{Terjoxin} proved the
generalized Jackson's inequalities for the cases of a bounded group
 and $s$-regular group. Remind that the group $\{U(t)\}_{t\in\R}$
is called $s$-regular if the resolvent of its generator $A$
satisfies the condition $\exists
\theta\in\R:\quad \|R_\lambda(e^{i\theta}A^s)\|\le\frac{C}{{\mathrm
Im}\lambda}$.

G. V. Radzievsky studied the direct and inverse theorems
\cite{Radzievsky1997, Radzievsky1998}, using the notion of
$K$-functional instead of module of continuity, but it should be
noted that the $K$-functional has two-sided estimates with regard to
module of continuity at least for bounded $C_0$-groups.

In the papers \cite{MGorbShilinst_ExpA,MGorb_OperAppr} and \cite{Gorb_Gr_Torba} the authors
investigated the case of a group of unitary operators in a Hilbert space and established
Jackson-type inequalities in Hilbert spaces and their rigs. These inequalities are used to estimate
the rate of convergence to zero
of the best approximation of both finite and infinite smoothness vectors for the operator $A$
 by exponential type entire vectors.

We consider the $C_0$-groups, generated by the so-called \emph{non-quasianalytic operators}
\cite{LubMatsaev}, i.e. the groups satisfying
\begin{equation}\label{NeKvaziAnalit}
\int_{-\infty}^{\infty}\frac{\ln\left\Vert U(t)\right\Vert
}{1+t^{2}}dt<\infty.
\end{equation}

As was shown in \cite{MGorbShilinst_ExpA}, the set of exponential type entire vectors
for the non-quasianalytic operator  $A$ is dense in
$\XX$, so the problem of approximation by exponential type entire vectors is
correct. On the other hand, it was shown in \cite{Gorbachuk_NeobhidnistNekvazianal} that condition
(\ref{NeKvaziAnalit}) is close to the necessary one, so in the case when
(\ref{NeKvaziAnalit}) doesn't hold, the class of entire vectors isn't necessary dense in $\XX$, and the
corresponding approximation problem loses its meaning.

The purpose of this work is to obtain Jackson-type inequalities in the case where
a vector of a Banach space is approximated by exponential type entire vectors for a non-quasianalytic
operator, and, in particular, Jackson-type inequalities in various classical function spaces.

\setcounter{equation}{0}
\section{Preliminaries}
Let $A$ be a closed linear operator with dense domain of definition
$\D(A)$ in the Banach space $(\XX,\left\Vert \cdot\right\Vert )$
over the field of complex numbers.

Let $C^{\infty}(A)$ denotes the set of all infinitely differentiable
vectors of the operator $A$, i.e.
\begin{equation*}
C^{\infty}(A)=\bigcap_{n\in\N_{0}}\D(A^{n}),\quad\N_{0}=\N\cup\{0\}.
\end{equation*}
For a number $\alpha>0$ we set
\begin{equation*}
\E^{\alpha}(A)=\left\{ x\in C^{\infty}(A)\,|\,\exists
c=c(x)>0\,\,\forall k\in\N_{0}\,\left\Vert A^{k}x\right\Vert \leq
c\alpha^{k}\right\}.
\end{equation*}
The set $\E^{\alpha}(A)$ is a Banach space with respect to the norm
\begin{equation*}
\left\Vert x\right\Vert
_{\E^{\alpha}(A)}=\sup_{n\in\N_{0}}\frac{\left\Vert
A^{n}x\right\Vert }{\alpha^{n}}\,.
\end{equation*}
Then $\E(A)=\bigcup_{\alpha>0}\E^{\alpha}(A)$ is a linear locally
convex space with respect to the topology of the inductive limit of
the Banach spaces $\E^{\alpha}(A)$:
\begin{equation*}
\E(A)=\limind_{\alpha\rightarrow\infty}\E^{\alpha}(A).
\end{equation*}
Elements of the space $\E(A)$ are called exponential type entire vectors of the operator $A$.
The type $\sigma(x,A)$ of a
vector $x\in\E(A)$ is defined as the number
\begin{equation*}
\sigma(x,A)=\inf\left\{ \alpha>0\,:\, x\in\E^{\alpha}(A)\right\}
=\limsup_{n\rightarrow\infty}\left\Vert A^{n}x\right\Vert
^{\frac{1}{n}}.
\end{equation*}

\begin{example}\label{AinLpDef}
Let $\XX$ is one of the $L_{p}(2\pi)$ ($1\leq p<\infty$) spaces of
integrable in $p$-th degree over $[0,2\pi]$, $2\pi$-periodical
functions or the space $C(2\pi)$ of continuous $2\pi$-periodical
functions (the norm in $\XX$ is defined in a standard way), and let
 $A$ is the differentiation operator in the space $\XX$ ($\D(A)=\{
x\in\XX\cap AC(\R)\,:\, x'\in\XX\}$; $(Ax)(t)=\frac{dx}{dt}$, where
$AC(\R)$ denotes the space of absolutely continuous functions over
$\R$). It can be proved that in such case the space $\E(A)$
coincides with the space of all trigonometric polynomials, and for
$y\in\E(A)$~ $\sigma(y,A)=\deg(y)$, where $\deg(y)$ is the degree of
the trigonometric polynomial $y$.
\end{example}

In what follows, we always assume that the operator $A$ is the
generator of the group of linear continuous operators $\{ U(t)\,:\,
t\in\R\}$ of class $C_{0}$ on $\XX$. We recall that belonging of
the group to the $C_{0}$ class means that for every $x\in\XX$ the
vector-function $U(t)x$ is continuous on $\R$ with respect to the
norm of the space $\XX$.

For $t\in\R_{+}$, we set
\begin{equation*}
M_{U}(t):=\sup_{\tau\in\mathbb{R},\,|\tau|\leq t}\left\Vert
U(\tau)\right\Vert.
\end{equation*}
The estimation $\|U(t)\|\le Me^{\omega t}$ for some $M,\omega\in\R$
implies $M_{U}(t)<\infty\,(\forall t\in\R_{+})$. It is easy to
see that the function $M_{U}(\cdot)$ has the following properties:
\begin{itemize}
\item[1)] $M_{U}(t)\geq1,$ $t\in\R_{+}$;
\item[2)] $M_{U}(\cdot)$ is monotonically non-decreasing on $\R_{+}$;
\item[3)] $M_{U}(t_{1}+t_{2})\leq M_{U}(t_{1})M_{U}(t_{2})$,
$t_{1},t_{2}\in\R_{+}$.
\end{itemize}

According to \cite{Kyptsov}, for $x\in\XX$, $t\in\R_{+}$ and
$k\in\N$ we set
\begin{gather}
\label{w_{d}ef} \omega_{k}(t,x,A)=\sup_{0\leq\tau\leq t}\left\Vert
\Delta_{\tau}^{k}x\right\Vert,\qquad\text{where}\\
 \label{Delta_{d}ef}
\Delta_{h}^{k}=(U(h)-\I)^{k}=\sum_{j=0}^{k}(-1)^{k-j}{{j \choose
k}}U(jh),\quad k\in\N_{0},\, h\in\R\quad (\Delta_{h}^{0}\equiv 1).
\end{gather}
Moreover, let
\begin{equation}
\widetilde{\omega}_{k}(t,x,A)=\sup_{|\tau|\leq t}\left\Vert
\Delta_{\tau}^{k}x\right\Vert.\label{w_sym_{d}ef}
\end{equation}

\begin{rmk}\label{RmkOmegaWhenU(t)Bounded}
It is easy to see that in the case of the bounded group $\{ U(t)\}$
($\left\Vert U(t)\right\Vert \leq M,\, t\in\R$) the quantities
 $\omega_{k}(t,x,A)$ and $\widetilde{\omega}_{k}(t,x,A)$
are equivalent within constant factor
($\omega_{k}(t,x,A)\leq\widetilde{\omega}_{k}(t,x,A)\leq
M\,\omega_{k}(t,x,A)$), and in the case of isometric group
($\left\Vert U(t)\right\Vert \equiv1,\, t\in\mathbb{R}$) these
quantities coincide.
\end{rmk}

It is immediate from the definition of
$\widetilde{\omega}_{k}(t,x,A)$ that for $k\in\N$:
\begin{itemize}
\item[1)] $\widetilde{\omega}_{k}(0,x,A)=0$;
\item[2)] for fixed $x$ the function $\widetilde{\omega}_{k}(t,x,A)$
is non-decreasing and is continuous by the variable $t$ on
$\mathbb{R}_{+}$;
\item[3)] $\widetilde{\omega}_{k}(nt,x,A)\leq\big(1+(n-1)M_{U}((n-1)t)\big)^{k}\widetilde{\omega}_{k}(t,x,A)$
~ ($n\in\N,\, t>0$);
\item[4)] $\widetilde{\omega}_{k}(\mu
t,x,A)\leq\big(1+\mu M_{U}(\mu
t)\big)^{k}\widetilde{\omega}_{k}(t,x,A)$ ~ ($\mu,t>0$);
\item[5)] for fixed $t\in\R_{+}$ the function
$\widetilde{\omega}_{k}(t,x,A)$ is continuous in $x$.
\end{itemize}

For arbitrary $x\in\XX$ we set, according to
\cite{Gorb_Gr_Torba,MGorb_OperAppr},
\begin{equation*}
\mathcal{E}_{r}(x,A)=\inf_{y\in\E(A)\,:\,\sigma(y,A)\leq
r}\left\Vert x-y\right\Vert ,\quad r>0,
\end{equation*}
i.e. $\mathcal{E}_{r}(x,A)$ is the best approximation of the element
 $x$ by exponential type entire vectors $y$ of the operator $A$ for which $\sigma(y,A)\leq r$.
For fixed $x$~ $\mathcal{E}_{r}(x,A)$ does not increase and
$\mathcal{E}_{r}(x,A)\rightarrow 0,\ r\rightarrow\infty$ for every
$x\in\XX$ if and only if the set $\E(A)$ of exponential type entire vectors is dense in $\XX$. Particularly, as indicated
above, the set $\E(A)$ is dense in $\XX$ if the group $\{ U(t)\,:\,
t\in\R\}$ belongs to non-quasianalytic class.

\setcounter{equation}{0}
\section{Abstract Jackson's inequality in a Banach space}
\begin{theorem}\label{Th1}
Suppose that $\{ U(t)\,:\, t\in\mathbb{R}\}$ satisfies condition
 (\ref{NeKvaziAnalit}). Then $\forall k\in\N$ there exists a constant $\mathbf{m}_{k}=\mathbf{m}_{k}(A)>0$,
 such that $\forall
x\in\XX$ the following inequality holds:
\begin{equation}\label{eq:2}
\mathcal{E}_{r}(x,A)\leq\mathbf{m}_{k}\cdot\tilde\omega_{k}\left(\frac{1}{r},x,A\right),\quad
r\ge 1.
\end{equation}
\end{theorem}

\begin{rmk}\label{Th1rmk}
If, additionally, the group $\left\{ U(t)\right\} $ is bounded
($M_{U}(t)\leq\widetilde{M}<\infty,\,\, t\in\mathbb{R}$), then the assumption $r\ge 1$ can be changed to $r>0$.
\end{rmk}

Integral kernels, constructed in \cite{Yadro_Marchenko}, will be
used in the proving of the theorem. Moreover, we need additional
properties of these kernels, lacking in \cite{Yadro_Marchenko}. The
following lemma shows how these kernels are constructed and
continues the investigation of their properties.

In what follows we denote as $\QQ$ the class of functions
 $\alpha:\R\mapsto\R$, satisfying the following conditions:

{\em
\begin{itemize}
\item [I)] $\alpha(\cdot)$ is measurable and bounded on any segment
$[-T,T]\subset\R$.
\item [II)] $\alpha(t)>0,\,\, t\in\R$.
\item [III)] $\alpha(t_{1}+t_{2})\leq\alpha(t_{1})\alpha(t_{2}),\quad t_{1},t_{2}\in\R$.
\item [IV)]
$\int_{-\infty}^{\infty}\frac{\left|\ln(\alpha(t))\right|}{1+t^{2}}dt<\infty$.
\end{itemize}
}
\begin{lemma}\label{Lemma_Yadro_1}
Let $\alpha\in\QQ$. Then there exists such entire function
$\K_{\alpha}:\CC\mapsto\CC$ that
\begin{itemize}
\item [1)] $\K_{\alpha}(t)\geq0,\,\, t\in\R$;
\item [2)] $\int_{-\infty}^{\infty}\K_{\alpha}(t)\,dt=1$; %\newcommand{\Imm}{{\bf Im}}
\item [3)] $\forall r>0\ \exists c_{r}=c_{r}(\alpha)>0\ \forall
z\in\CC\quad |\K_{\alpha}(rz)|\leq c_{r}\frac{e^{r|\Imm
z|}}{\alpha(|z|)}$
\end{itemize}
\end{lemma}
\begin{proof}
Without lost of generality we may assume that the function
$\alpha(t)$ satisfies additional conditions:
{\em
\begin{itemize}
\item [V)] $\alpha(t)\ge 1$, $t\in\R$;
 \footnote{As shown in \cite{LubMatsaev}, for non-quasianalytic groups
 the condition $\|U(t)\|\ge 1$ always holds, therefore in this paper
 the condition V) automatically takes place.}
\item [VI)] $\alpha(t)$ is even on $\R$ and is monotonically increasing on
$\R_{+}$;
\item[VII)] $\left\Vert \alpha^{-1}\right\Vert _{L_{1}(\R)}=\int_{-\infty}^{\infty}|\alpha^{-1}(t)|dt<\infty$.
\end{itemize}
}

It is easy to verify that assumptions V),VII) and condition that
the function $\alpha(t)$ is even in VI) don't confine the general
case if one examined the function
$\alpha_{1}(t)=\widetilde{\alpha}(t)\widetilde{\alpha}(-t)$, where
$\widetilde{\alpha}(t)=(1+\alpha(t))(1+t^{2})$. In \cite[theorems 1
and 2]{Marchenko2} it has been proved that the monotony condition on
$\alpha(t)$ in VI) doesn't confine the general case too.

It follows from VII) that
\begin{equation}\label{AlphaToInfty}
\alpha(t)\to\infty,\ t\to\infty.
\end{equation}

Let $\beta(t)=\ln\alpha(t),\,\, t\in\R$. Conditions III)-VII) and
(\ref{AlphaToInfty}) lead to conclusion that
\begin{equation*}
\beta(t)>0, \qquad \beta(-t)=\beta(t),\qquad \beta(t)\to\infty,\,\,
t\to\infty;
\end{equation*}
\begin{align}
 \beta(t_{1}+t_{2})& \leq\beta(t_{1})+\beta(t_{2}),\,\, t_{1},t_{2}\in\R\label{NapivadBeta}\\
 \int_{1}^{\infty}\frac{\beta(t)}{t^{2}}dt & <\infty\label{IntBeta}
\end{align}
Because of (\ref{NapivadBeta}) there exists limit
$\lim_{t\to\infty}\frac{\beta(t)}{t}$. And, by virtue of
(\ref{IntBeta}):
\begin{equation}
\lim_{t\to\infty}\frac{\beta(t)}{t}=0.\label{LimBeta}
\end{equation}
Also, using (\ref{IntBeta}) it is easy to check that
\begin{equation}\label{SumBeta}
\sum_{k=1}^{\infty}\frac{\beta(k)}{k^{2}}<\infty,
\end{equation}
moreover, all terms of the series (\ref{SumBeta}) are positive. From
the convergence of series (\ref{SumBeta}) follows the existence
of such sequence  $\{ Q_{n}\}_{n=1}^{\infty}\subset\R$ that
$Q_{n}>1,\quad Q_{n}\to\infty,\ n\to\infty$ and
\begin{equation}
\sum_{k=1}^{\infty}\frac{\beta(k)}{k^{2}}Q_{k}=S<\infty.\label{SumBetaQ}
 \end{equation}

We set
\begin{equation*}
a_{k}:=\frac{\beta(k)Q_{k}}{S\, k^{2}},\quad k\in\N.
\end{equation*}
The definition of $a_k$ and (\ref{SumBetaQ}) result in equality
\begin{equation}\label{S_a_k}
\sum_{k=1}^{\infty}a_{k}=1.
\end{equation}
We construct the sequence of functions, which, obviously, are entire
for every $n\in\N$:
\begin{equation*}
f_{n}(z):=\prod_{k=1}^{n}P_{k}(z),\quad\textrm{where}\,\,
P_{k}(z)=\left(\frac{\sin\frac{a_{k}z}{2}}{\frac{a_{k}z}{2}}\right)^{2},\,\,
z\in\CC,\, n\in\N.
\end{equation*}

Similarly to the proof of the Denjoy-Carleman theorem
\cite[p.378]{Rudin1970} it can be concluded that the sequence of
(entire) functions $f_{n}(z)$ converges uniformly to the function
\begin{equation*}
f(z)=\prod_{k=1}^{\infty}\left(\frac{\sin\frac{a_{k}z}{2}}{\frac{a_{k}z}{2}}\right)^{2},\qquad
z\in\CC
\end{equation*}
in every disk $\{ z\in\CC\,|\,|z|\leq R\}$. Thus, by Weierstrass
theorem, the function $f(z)$ is entire.

Using the inequality $|\sin z|\leq\min(1,|z|)e^{|\Imm z|}$,
$z\in\CC$ and taking (\ref{S_a_k}) into account, when $z\in\CC$ and
$r>0$, we receive
\begin{multline*}
\left|f(rz)\right|=\prod_{k=1}^{\infty}\left|\frac{\sin\frac{a_{k}rz}{2}}{\frac{a_{k}rz}{2}}\right|^{2}\leq\prod_{k=1}^{\infty}\left(\frac{2}{a_{k}r|z|}\min\left(1,\frac{a_{k}r|z|}{2}\right)e^{\frac{1}{2}a_{k}r|\Imm z|}\right)^{2}=\\
=e^{r|\Imm
z|}\prod_{k=1}^{\infty}{\min}^{2}\left(1,\frac{2}{a_{k}r|z|}\right)\le
e^{r|\Imm z|}\prod_{k=1}^N{\min}^2\left(1, \frac 2{a_{k}r|z|}\right)
\end{multline*}
 for every $N\in\N$. Using the inequality
 $\min(1,a)\cdot\min(1,b)\le\min(1,ab)$, we get:
\begin{multline}\label{f(rz)_1}
\left|f(rz)\right|\leq e^{r|\Imm
z|}{\min}^{2}\left(1,\prod_{k=1}^{N}\frac{2}{a_{k}r|z|}\right)=e^{r|\Imm z|}{\min}^{2}\Bigg(1,\frac{2^{N}}{\left(\prod_{k=1}^{N}\frac{\beta(k)Q_{k}}{S\, k^{2}}\right)(r|z|)^{N}}\Bigg)=\\
=e^{r|\Imm
z|}{\min}^{2}\bigg(1,\frac{2^{N}N!}{\frac{\beta(1)}{1}\cdots\frac{\beta(N)}{N}\,\left(\frac{r}{S}\right)^{N}|z|^{N}Q_{1}\cdots
Q_{N}}\bigg).
 \end{multline}
Because of the condition $Q_{n}\to\infty$, $n\to\infty$ there exists
such number $n(r)\in\N$ that:
\begin{equation}\label{n(r)}
\forall n>n(r)\quad Q_{n}\geq\frac{4\sqrt{e}S}{r}.
\end{equation}
It follows from (\ref{LimBeta}) that there is
$T_{0}\in(0,\infty)$ such that:
\begin{equation}\label{Bt/t_leq_1}
\forall\, t>T_{0}\quad\frac{\beta(t)}{t}\leq1.
\end{equation}
In \cite{Yadro_Marchenko} the following statement was proved:
\begin{equation}\label{QuaziMonot_Bt/t}
\forall\, t_{1},t_{2}\in\R_{+}\quad t_{1}\leq
t_{2}\,\Rightarrow\,\frac{\beta(t_{1})}{t_{1}}\geq\frac{1}{2}\frac{\beta(t_{2})}{t_{2}}.
\end{equation}

Let $z\in\CC$ and $|z|\geq\max\big(\beta^{[-1]}(n(r)),T_{0}\big),$
where $\beta^{[-1]}$ is the inverse function of the function $\beta$
on $[0,\infty)$ (the inverse function $\beta^{[-1]}$ exists due to
monotony of $\beta$ on $[0,\infty)$). We substitute as $N$ in
(\ref{f(rz)_1}) $N:=[\beta(|z|)]$, where $[\cdot]$ denotes the
integer part of a number. Then for $k\in\{1,\dots,N\}$, in
accordance with (\ref{Bt/t_leq_1}) and (\ref{QuaziMonot_Bt/t}), we
obtain $k\leq N\leq\beta(|z|)\leq|z|$ and
\begin{equation}
 \frac{\beta(k)}{k}\ge\frac{1}{2}\frac{\beta(|z|)}{|z|}.
 \label{Bk/k}
\end{equation}
Using (\ref{f(rz)_1}),(\ref{n(r)}),(\ref{Bk/k}), we find
\begin{multline*}
  |f(rz)|\leq e^{r|\Imm z|}\Bigg(\frac{2^{N}N!}{\big(\frac{1}{2}\frac{\beta(|z|)}{|z|}\big)^{N}\,\left(\frac{r}{S}\right)^{N}|z|^{N}Q_{1}\cdots Q_{N}}\Bigg)^{2}\leq  \displaybreak[1]\\
  \leq e^{r|\Imm z|}\Bigg(\frac{2^{N}N!}{\big(\frac{1}{2}\frac{N}{|z|}\big)^{N}\,\left(\frac{r}{S}\right)^{N}|z|^{N}Q_{1}\cdots Q_{N}}\Bigg)^{2}=e^{r|\Imm z|}\left(\frac{2^{2N}N!}{N^{N}\,\left(\frac{r}{S}\right)^{N}Q_{1}\cdots Q_{N}}\right)^{2}\leq \displaybreak[1]\\
  \leq e^{r|\Imm z|}\left(\frac{2^{2N}}{\left(\frac{r}{S}\right)^{N}Q_{1}\cdots Q_{N}}\right)^{2}=e^{r|\Imm z|}\left(\frac{\left(\frac{4S}{r}\right)^{N}}{Q_{1}\cdots
  Q_{N}}\right)^{2}.
\end{multline*}
Since $Q_n\ge 1$, the last inequality leads to
\begin{multline}\label{f(rz)_2}
  |f(rz)|\leq e^{r|\Imm z|}\Bigg(\frac{\left(\frac{4S}{r}\right)^{N}}{\big(\frac{4\sqrt{e}S}{r}\big)^{N-n(r)}}\Bigg)^{2}=e^{r|\Imm z|}\left(\frac{4\sqrt{e}S}{r}\right)^{2n(r)}e^{-[\beta(|z|)]}\leq\\
  \le e^{r|\Imm z|}\left(\frac{4\sqrt{e}S}{r}\right)^{2n(r)}e^{-(\beta(|z|)-1)}= C_{r}^{(1)}\frac{e^{r|\Imm z|}}{\alpha(|z|)},
 \end{multline}
where $C_{r}^{(1)}=e\left(\frac{4\sqrt{e}S}{r}\right)^{2n(r)}$. When
$z\in\CC$ and $|z|<\max\big(\beta^{[-1]}(n(r)),T_{0}\big)$, using
(\ref{f(rz)_1}), we get
\begin{equation}\label{f(rz)_3}
 |f(rz)|\leq e^{r|\Imm z|}=e^{r|\Imm
z|}\frac{\alpha(|z|)}{\alpha(|z|)}\leq e^{r|\Imm
z|}\frac{C_{r}^{(2)}}{\alpha(|z|)},
\end{equation}
where $C_{r}^{(2)}=\alpha(\max(\beta^{[-1]}(n(r)),T_{0}))$. It
follows from (\ref{f(rz)_2}), (\ref{f(rz)_3}) that
\begin{equation}\label{f(rz)}
|f(rz)|\leq e^{r|\Imm z|}\frac{C_{r}^{(0)}}{\alpha(|z|)},\quad
z\in\CC,\qquad\textrm{where}\quad
C_{r}^{(0)}=\max(C_{r}^{(1)},C_{r}^{(2)}).
\end{equation}

Inequality (\ref{f(rz)}) and Condition VII) imply that
$\left\Vert f\right\Vert _{L_{1}(\R)}<\infty$. Thus it is enough to
set $\K_{\alpha}(z):=\frac{1}{\left\Vert f\right\Vert
_{L_{1}(\R)}}f(z)$, $z\in\CC$ and use (\ref{f(rz)}) to finish the
proof.
\end{proof}

Let $\alpha\in\QQ$, and ~$\K_{\alpha}:\CC\mapsto\CC$ is the function
constructed by the function $\alpha$ in lemma \ref{Lemma_Yadro_1}.
We set
\begin{equation*}
\K_{\alpha,r}(z):=r\K_{\alpha}(rz),\quad z\in\CC,\,\,
r\in(0,\infty).
\end{equation*}
The lemma \ref{Lemma_Yadro_1} ensures us that the function
 $\K_{\alpha,r}$ has the following properties:
\begin{itemize}
\item [1)] $\K_{\alpha,r}(t)\geq0,\quad t\in\R$;
\item [2)] $\int_{-\infty}^{\infty}\K_{\alpha,r}(t)\,dt=1$;
\item [3)] $\forall z\in\CC\,\,\,|\K_{\alpha,r}(z)|\leq rc_{r}\frac{e^{r|\Imm z|}}{\alpha(|z|)};\quad r>0$.
\end{itemize}
\begin{lemma}\label{Lemma_Yadro_2}
$\forall r\in(0,\infty)$ there exists constant
$\cc_{r}=\cc_{r}(\alpha)>0$, such that $\forall n\in\N$ the
following inequality holds:
\begin{equation*}
|\K_{\alpha,r}^{(n)}(t)|\leq\cc_{r}\frac{\sqrt{2\pi
n}\,\alpha\left(\frac{n}{r}\right)}{\alpha(|t|)}r^{n},\qquad t\in\R
\end{equation*}
\end{lemma}
\begin{proof}
In what follows in this proof we assume $t\in\R$, $r\in(0,\infty)$,
$n\in\N$. Let
\begin{equation*}
\gamma_{n,r}(t):=\left\{
\zeta\in\CC\,:\,|\zeta-t|=\frac{n}{r}\right\}.
\end{equation*}

Using Cauchy's integral theorem and Stirling's approximation for
$n!$, we get
\begin{multline*}
 |\K_{\alpha,r}^{(n)}(t)|\leq\frac{n!}{2\pi}\,\oint_{\gamma_{n,r}(t)}\frac{|\K_{\alpha,r}(\xi)|}{|\xi-t|^{n+1}}|d\xi|=\frac{n!}{2\pi}\frac{r^{n+1}}{n^{n+1}}\,\oint_{\gamma_{n,r}(t)}|\K_{\alpha,r}(\xi)||d\xi|\leq\\
 \leq\frac{c^{(!)}r^{n+1}}{\sqrt{2\pi n}}e^{-n}\oint_{\gamma_{n,r}(t)}|\K_{\alpha,r}(\xi)||d\xi|,\quad \text{where}\ c^{(!)}=\sup_{k\in\N}\,\frac{k!}{\sqrt{2\pi
 k}}\left(\frac ek\right)^k< e^{1/12}.
\end{multline*}
Using property 3) of the function $K_{\alpha,r}$, the condition
$t\in\R$ and conditions III), VI) of the function $\alpha$, one
can find from the last inequality
\begin{multline*}
|\K_{\alpha,r}^{(n)}(t)|\leq\frac{c^{(!)}r^{n+1}}{\sqrt{2\pi n}}e^{-n}rc_{r}\oint_{\gamma_{n,r}(t)}\frac{e^{r|\Imm\xi|}}{\alpha(|\xi|)}|d\xi|=\\
=\frac{c^{(!)}r^{n+1}}{\sqrt{2\pi n}}e^{-n}\frac{rc_{r}}{\alpha(|t|)}\oint_{\gamma_{n,r}(t)}\frac{e^{r|\Imm(\xi-t)|}\alpha(|(t-\xi)+\xi|)}{\alpha(|\xi|)}|d\xi|\leq\\
 \leq\frac{c^{(!)}r^{n+1}}{\sqrt{2\pi n}}e^{-n}\frac{rc_{r}}{\alpha(|t|)}\oint_{\gamma_{n,r}(t)}e^{r|\Imm(\xi-t)|}\alpha(|t-\xi|)|d\xi|\le\\
 \le\frac{c^{(!)}r^{n+1}}{\sqrt{2\pi
n}}e^{-n}\frac{rc_{r}}{\alpha(|t|)}\oint_{\gamma_{n,r}(t)}e^{n}\alpha\left(\frac{n}{r}\right)|d\xi|=\cc_{r}\frac{\sqrt{2\pi
n}\,\alpha\left(\frac{n}{r}\right)}{\alpha(|t|)}r^{n},
\end{multline*}
where $\cc_{r}=c^{(!)}rc_{r}$.
\end{proof}

\begin{rmk}\label{RmkLemma_Yadro_1}
If the function $\alpha(t)$ satisfies the conditions of lemma
\ref{Lemma_Yadro_1}, but, moreover, has the polynomial order of
growth at infinity, i.e. $\exists m\in\N_0,\,\exists M>0$:
\begin{equation}\label{Lemms_remark}
    \alpha(t)\le M(1+|t|)^{2m},\quad  t\in\R,
\end{equation}
another integral kernel may be used:
\begin{equation*}
    \tilde K_\alpha(z)=\frac 1{\textrm K_m}\left(\frac{\sin\frac
    z{2m}}{\frac z{2m}}\right)^{2m},\qquad {\textrm
    K_m}=\int_{-\infty}^\infty\left(\frac{\sin\frac
    x{2m}}{\frac x{2m}}\right)^{2m}dx.
\end{equation*}

In much the same way to the proving of the lemmas
\ref{Lemma_Yadro_1} and \ref{Lemma_Yadro_2} one can show that
\begin{equation*}
    \big|\tilde K_{\alpha}(rz)\big|\le \tilde C_r\frac{e^{r|\Imm
    z|}}{\alpha(|z|)},\qquad\text{where}\ \tilde C_r=\frac M{\textrm
    K_m}\left(1+\frac {2m}r\right)^{2m},
\end{equation*}
and
\begin{equation*}
    \big|\tilde K_{\alpha,r}^{(n)}(t)\big|\le \tilde c_r\frac{\sqrt{2\pi n}\,\alpha(\frac
    nr)}{\alpha(|t|)}r^n,\qquad\text{where}\ \tilde c_r=c^{(!)}r\tilde
    C_r,
\end{equation*}
that is to say, defined in such a way integral kernel satisfies
lemmas \ref{Lemma_Yadro_1} and \ref{Lemma_Yadro_2}.
\end{rmk}

\begin{proof}[Proof of theorem \ref{Th1}]
Let the group $\{ U(t):t\in\R\}$ satisfies (\ref{NeKvaziAnalit}). Then
it follows from \cite[theorems 1 and 2]{Marchenko2} that
\begin{equation}\label{L0L}
\int_{-\infty}^{\infty}\frac{\ln\left(M_{U}(|t|)\right)}{1+t^{2}}dt<\infty.
\end{equation}
 We fix arbitrary $k\in\N$ and set
\begin{equation*}
\alpha(t):=\big(M_{U}(|t|)\big)^{k}(1+|t|)^{k+2},\qquad t\in\R.
\end{equation*}
The function $\alpha$ is, obviously, even on $\R$. Condition
(\ref{L0L}) and the properties of the function $M_{U}(\cdot)$ imply
 $\alpha\in\QQ$, and, moreover,
\begin{equation}\label{UmovaNaAlpha}
\int_{-\infty}^{\infty}\frac{\big((1+|t|)M_{U}(|t|)\big)^{k}}{\alpha(t)}dt=\int_{-\infty}^{\infty}\frac{dt}{(1+|t|)^{2}}=2.
\end{equation}
Using lemma \ref{Lemma_Yadro_1} (or remark
\ref{RmkLemma_Yadro_1} if $\alpha(t)\le M(1+|t|)^m$) for the
function $\alpha(t)$, we construct the family of kernels
$K_{\alpha,r}$.

In what follows, we assume $x\in\XX,\ r\in (0,\infty)$ and
$n\in\{1,\ldots,k\}$. We define
\begin{equation*}
x_{r,n}:=\int_{-\infty}^{\infty}\K_{\alpha,r}(t)U(nt)x\,dt.
\end{equation*}
Let $\nu\in\N_{0}$. Let's prove that $x_{r,n}\in
C^{\infty}(A)=\bigcap_{\nu\in\N_{0}}\D(A^{\nu})$~ and
\begin{equation}
A^{\nu}x_{r,n}=\frac{(-1)^{\nu}}{n^{\nu}}\int_{-\infty}^{\infty}\K_{\alpha,r}^{(\nu)}(t)U(nt)x\,dt.\label{x(rn)inC(infty)}
\end{equation}
It follows from the property 3) of the function $\K_{\alpha,r}$ and
from lemma \ref{Lemma_Yadro_2} that there exists such constant
$\widetilde{C}(\nu,r)>0$ that
$\K_{\alpha,r}^{(\nu)}(t)\leq\frac{\widetilde{C}(\nu,r)}{\alpha(t)}$,
$t\in\R$. Thus, using (\ref{UmovaNaAlpha}), we get
\begin{multline}\label{IntK(nu)U}
\int_{-\infty}^{\infty}\left\Vert \K_{\alpha,r}^{(\nu)}(t)U(nt)x\right\Vert dt\leq\int_{-\infty}^{\infty}\frac{\widetilde{C}(\nu,r)}{\alpha(t)}\left\Vert U(t)\right\Vert ^{n}\left\Vert x\right\Vert dt\leq\\
\le\widetilde{C}(\nu,r)\left\Vert x\right\Vert
\int_{-\infty}^{\infty}\frac{M_{U}(|t|)^{k}}{\alpha(t)}dt\leq
2\widetilde{C}(\nu,r)\left\Vert x\right\Vert <\infty.
\end{multline}
Therefore the integral
$\int_{-\infty}^{\infty}\K_{\alpha,r}^{(\nu)}(t)U(nt)x\,dt$
converges. We define
\begin{equation*}
x_{r,n}^{(\nu)}=\frac{(-1)^{\nu}}{n^{\nu}}\int_{-\infty}^{\infty}\K_{\alpha,r}^{(\nu)}(t)U(nt)x\,dt.
\end{equation*}
Then, using closedness of the operator $A$ and integration by parts,
one can find for $x\in\D(A)$ that $x_{r,n}^{(\nu)}\in\D(A)$ and
\begin{multline}\label{Ax(nu)D(A)}
 Ax_{r,n}^{(\nu)}=\frac{(-1)^{\nu}}{n^{\nu}}\int_{-\infty}^{\infty}\K_{\alpha,r}^{(\nu)}(t)U(nt)Ax\,dt=\frac{(-1)^{\nu}}{n^{\nu}}\frac{1}{n}\int_{-\infty}^{\infty}\K_{\alpha,r}^{(\nu)}(t)(U(nt)x)'dt= \\
 =-\frac{(-1)^{\nu}}{n^{\nu}}\frac{1}{n}\int_{-\infty}^{\infty}\K_{\alpha,r}^{(\nu+1)}(t)U(nt)x\,dt=x_{r,n}^{(\nu+1)}.
\end{multline}
Let $x$ is an arbitrary element of the space $\XX$. Then there
exists the sequence $\{ x_{m}\}_{m=1}^{\infty}\subset\D(A)$ such
that $\left\Vert x_{m}-x\right\Vert \to0,\quad m\to\infty$.
Consequently, using inequality (\ref{IntK(nu)U}) and
relation (\ref{Ax(nu)D(A)}), one can get
\begin{gather*}
\left\Vert (x_{m})_{r,n}^{(\nu)}-x_{r,n}^{(\nu)}\right\Vert
\leq\frac{1}{n^{\nu}}\int_{-\infty}^{\infty}\left\Vert
\K_{\alpha,r}^{(\nu)}(t)U(nt)(x_{m}-x)\right\Vert dt
\leq\frac{2\widetilde{C}(\nu,r)}{n^{\nu}}\left\Vert
x_{m}-x\right\Vert \to 0;\\
 \left\Vert A(x_{m})_{r,n}^{(\nu)}-x_{r,n}^{(\nu+1)}\right\Vert =\left\Vert
(x_{m})_{r,n}^{(\nu+1)}-x_{r,n}^{(\nu+1)}\right\Vert \to0,\quad
m\to\infty.
\end{gather*}
Hence, taking into account closedness of the operator $A$, we
have:
\begin{equation}\label{Ax(nu)}
x_{r,n}^{(\nu)}\in\D(A),\quad Ax_{r,n}^{(\nu)}=x_{r,n}^{(\nu+1)}.
\end{equation}
One can get (\ref{x(rn)inC(infty)}) from (\ref{Ax(nu)}) by
induction.

Using relation (\ref{x(rn)inC(infty)}) and lemma
\ref{Lemma_Yadro_2}, one can find:
\begin{multline}\label{NormA(nu)x}
\left\Vert A^{\nu}x_{r,n}\right\Vert \leq\frac{\left\Vert x\right\Vert }{n^{\nu}}\int_{-\infty}^{\infty}\left|\K_{\alpha,r}^{(\nu)}(t)\right|\left\Vert U(nt)\right\Vert dt\leq\\
\leq\frac{\left\Vert x\right\Vert
}{n^{\nu}}\int_{-\infty}^{\infty}\cc_{r}\frac{\sqrt{2\pi\nu}\,\alpha\left(\frac{\nu}{r}\right)}{\alpha(|t|)}r^{\nu}\left\Vert
U(t)\right\Vert ^{n}dt\leq\\
\leq\cc_{r}\left\Vert x\right\Vert
\sqrt{2\pi\nu}\,\alpha\left(\frac{\nu}{r}\right)\bigg(\int_{-\infty}^{\infty}\frac{\left\Vert
U(t)\right\Vert
^{n}}{\alpha(t)}dt\bigg)\left(\frac{r}{n}\right)^{\nu},
\end{multline}
 where, accordingly to (\ref{UmovaNaAlpha}) and due to $n\le k$, $\int_{-\infty}^{\infty}\frac{\left\Vert U(t)\right\Vert ^{n}}{\alpha(t)}dt\leq\int_{-\infty}^{\infty}\frac{\left\Vert U(t)\right\Vert ^{k}}{\alpha(t)}dt\leq2<\infty$.
Since $\beta(t)=\ln(\alpha(t))$, $t\in\R$, as was mentioned in the
proof of lemma \ref{Lemma_Yadro_1},
$\lim_{\tau\to\infty}\frac{\beta(\tau)}{\tau}=0$ (cf.
(\ref{LimBeta})). Thus
\begin{equation*}
\lim_{\nu\to\infty}\left(\alpha\left(\frac{\nu}{r}\right)\right)^{1/\nu}=\lim_{\nu\to\infty}e^{\frac{1}{r}\left(\frac{r}{\nu}\beta\left(\frac{\nu}{r}\right)\right)}=e^{\frac{1}{r}\cdot0}=1.
\end{equation*}
Therefore from relation (\ref{NormA(nu)x}) one can get:
\begin{equation*}
\limsup_{\nu\to\infty}\big(\left\Vert A^{\nu}x_{r,n}\right\Vert
\big)^{1/\nu}\leq\frac{r}{n}.
\end{equation*}
The last inequality brings us to the conclusion that
\begin{equation}\label{x(rn)inE(A)}
x_{r,n}\in\E(A)\quad\textrm{and}\quad\sigma(x_{r,n},A)\leq\frac{r}{n}.
 \end{equation}

For arbitrary  $x\in\XX$ we define
\begin{multline}\label{xx_rk(def)}
 \xx_{r,k}:=\int_{-\infty}^{\infty}\K_{\alpha,r}(t)(x+(-1)^{k-1}(U(t)-\I)^{k}x)dt= \\
 =\int_{-\infty}^{\infty}\K_{\alpha,r}(t)\sum_{n=1}^{k}(-1)^{n+1}{{k \choose
n}}U(nt)x\,dt
\end{multline}
(the absolute convergence by the norm of $\XX$ of the integral in
the right part of (\ref{xx_rk(def)}) follows from inequality
(\ref{IntK(nu)U}), so the definition of the vector $\xx_{r,k}$ is
correct). Using definition (\ref{xx_rk(def)}) one can get:
\begin{equation*}
\xx_{r,k}=\sum_{n=1}^{k}(-1)^{n+1}{{k \choose
n}}\int_{-\infty}^{\infty}\K_{\alpha,r}(t)U(nt)xdt=\sum_{n=1}^{k}(-1)^{n+1}{{k
\choose n}}x_{r,n}.
\end{equation*}
Therefore, accordingly to (\ref{x(rn)inE(A)}),
\begin{equation*}
\xx_{r,k}\in\E(A)\quad\textrm{and}\quad\sigma(\xx_{r,k},A)\leq r.
\end{equation*}
Hence for an arbitrary  $x\in\XX$ we have:
\begin{equation*}
\mathcal{E}_{r}(x,A)=\inf_{y\in\E(A)\,:\,\sigma(y,A)\leq
r}\,\left\Vert x-y\right\Vert \leq\left\Vert x-\xx_{r,k}\right\Vert
\end{equation*}
Using (\ref{xx_rk(def)}), the property 2) of the kernel
$\K_{\alpha,r}$ and (\ref{w_sym_{d}ef}), the last inequality
implies:
\begin{multline*}
\mathcal{E}_{r}(x,A)\leq\left\Vert
\int_{-\infty}^{\infty}\K_{\alpha,r}(t)x
dt-\int_{-\infty}^{\infty}\K_{\alpha,r}(t)\big(x+(-1)^{k-1}(U(t)-\I)^{k}x\big)\,dt\right\Vert\le\\
\leq\int_{-\infty}^{\infty}\K_{\alpha,r}(t)\left\Vert
(U(t)-\I)^{k}x\right\Vert
dt\leq\int_{-\infty}^{\infty}\K_{\alpha,r}(t)\widetilde{\omega}_{k}(|t|,x,A)\,dt.
 \end{multline*}
So, in accordance with the property 4) of the function
$\widetilde{\omega}_{k}(|t|,x,A)$,
\begin{multline}\label{Er(1)}
 \mathcal{E}_{r}(x,A)\leq\int_{-\infty}^{\infty}\K_{\alpha,r}(t)\widetilde{\omega}_{k}\left(|rt|\frac{1}{r},x,A\right)dt\leq \\
 \leq\widetilde{\omega}_{k}\left(\frac{1}{r},x,A\right)\int_{-\infty}^{\infty}\big(1+|rt|M_{U}(|t|)\big)^{k}\K_{\alpha,r}(t)dt.
 \end{multline}
Taking into account properties of the function $M_{U}(\cdot)$,
the definition of $\K_{\alpha,r}$, lemma \ref{Lemma_Yadro_1} and
 equality (\ref{UmovaNaAlpha}), one can find for $r\ge 1$:
\begin{multline*}
\int_{-\infty}^{\infty}\big(1+|rt|M_{U}(|t|)\big)^{k}\K_{\alpha,r}(t)dt\leq\int_{-\infty}^{\infty}\big(1+|rt|M_{U}(rt)\big)^{k}r\K_{\alpha}(rt)dt\leq\\
\leq\int_{-\infty}^{\infty}\big((1+\tau)M_{U}(\tau)\big)^{k}\K_{\alpha}(\tau)d\tau\leq c_{1}\int
_{-\infty}^{\infty}\frac{\big((1+|\tau|)M_{U}(|\tau|)\big)^{k}}{\alpha(\tau)}d\tau=2c_{1}<\infty.
\end{multline*}
In accordance with (\ref{Er(1)}), inequality (\ref{eq:2}) holds
for all $r\in[1,\infty)$ with a constant $\mathbf{m}_{k}=2c_{1}$. It
should be noted that constant $\mathbf{m}_{k}$, indeed, depends on
$k$, because due to \ref{Lemma_Yadro_1}, the constant
$c_{1}=c_{1}(\alpha)$ depends on the function
$\alpha(t)=(M_{U}(|t|))^{k}(1+|t|)^{k+2}$.

Moreover, let the group $\{U(t)\}$ is bounded ($M_{U}(t)\leq \widetilde{M},\, \, t\in \mathbb{R}$,
$\widetilde{M}\ge 1$). Taking into account properties of the function $M_{U}(\cdot)$,
the definition of $\K_{\alpha,r}$, lemma \ref{Lemma_Yadro_1} and
equality (\ref{UmovaNaAlpha}), one can find for $r\in(0,\infty)$
\begin{multline*}
\int_{-\infty }^{\infty }(1+|rt|M_{U}(|t|))^{k}\K _{\alpha ,r}(t)dt\leq \int_{-\infty }^{\infty }(1+|rt|\widetilde{M}M_{U}(rt))^{k}r\K _{\alpha }(rt)dt\leq \\
\le\widetilde{M}^{k}\int_{-\infty }^{\infty }((1+\tau )M_{U}(\tau ))^{k}\K _{\alpha }(\tau )d\tau \leq 2\widetilde{M}^{k}c_{1}<\infty,
\end{multline*}
which proves remark \ref{Th1rmk} with the constant $\mathbf{m}_{k}=2\widetilde{M}^{k}c_{1}$.
\end{proof}

Theorem \ref{Th1} allows us to prove the analogue of the classic
Jackson's inequality for $m$ times differentiable functions:
\begin{nasl}\label{NaslDlyaD(A^m)}
Let $x\in\D(A^{m}),\,\, m\in\N_{0}$. Then $\forall k\in\N_0$
\begin{equation}
\mathcal{E}_{r}(x,A)\leq\mathbf{m}_{k+m}\frac{M_{U}\left(\frac{m}{r}\right)}{r^{m}}\widetilde{\omega}_{k}\left(\frac{1}{r},A^{m}x,A\right),\quad
r\ge 1,\label{eq:3}
\end{equation}
where the constants $\mathbf{m}_{n}$ ($n\in\N$) are the same as in
 theorem \ref{Th1}.
\end{nasl}
\begin{proof}
Let $x\in\D(A^{m})$ and $r\ge 1$. By theorem \ref{Th1},
\begin{equation*}
\mathcal{E}_{r}(x,A)\leq\mathbf{m}_{k+m}\cdot\widetilde{\omega}_{k+m}\left(\frac{1}{r},x,A\right).
\end{equation*}
Let $t\in\R$, $0\leq|t|\leq\frac{1}{r}$. Then, using properties of
the groups of the $C_{0}$ class and properties of the function
$M_{U}(t)$, one can get:
\begin{multline*}
 \left\Vert (U(t)-\I)^{k+m}x\right\Vert =\left\Vert (U(t)-\I)^{m}(U(t)-\I)^{k}x\right\Vert \leq\\
 \leq\int_{0}^{t}\cdots\int_{0}^{t}\left\Vert
U(\xi_{1}+\dots+\xi_{m})\right\Vert \left\Vert
(U(t)-\I)^{k}A^{m}x\right\Vert \, d\xi_{1}\dots d\xi_{m}\leq \\
\leq M_{U}(m|t|)\left\Vert (U(t)-\I)^{k}A^{m}x\right\Vert
t^{m}\leq\frac{M_{U}(\frac{m}{r})}{r^{m}}\widetilde{\omega}_{k}\left(\frac{1}{r},A^{m}x,A\right).
\end{multline*}
This implies $\widetilde{\omega}_{k+m}\left(\frac{1}{r},x,A\right)=
\sup_{|t|\leq\frac{1}{r}}\left\Vert (U(t)-\I)^{k+m}x\right\Vert
\leq\frac{M_{U}(\frac{m}{r})}{r^{m}}\widetilde{\omega}_{k}\left(\frac{1}{r},A^{m}x,A\right),$
which proves inequality (\ref{eq:3}).
\end{proof}

By setting in corollary \ref{NaslDlyaD(A^m)} $k=0$ and taking into account that
$\tilde\omega_0\left(\cdot,A^{m}x,A\right)\equiv \|A^{m}x\|$,
 one can conclude the following inequality:
\begin{nasl}\label{NaslDlyaD(A^m)2}
Let $x\in\D(A^{m}),\,\, m\in\N_{0}$. Then
\begin{equation}
\mathcal{E}_{r}(x,A)\leq\frac{\mathbf{m}_{m}}{r^{m}}\big(M_{U}(1/r)\big)^{m}\|A^{m}x\|\qquad
r\ge 1,\label{eq:4}
\end{equation}
where the constants $\mathbf{m}_{n}$ ($n\in\N$) are the same as in
theorem \ref{Th1}.
\end{nasl}

\setcounter{equation}{0}
\section{The examples of application of the abstract Jackson's inequality in particular spaces}
Lets consider several examples of application of theorem
 \ref{Th1} in particular spaces.

\subsection{Jackson's inequalities in $L_{p}(2\pi)$ and $C(2\pi)$}

\begin{example}\label{Example2}
Let the space $\XX$ and the operator $A$ are the same as in the
example \ref{AinLpDef}. Then for $x\in\XX$ the quantity
$\mathcal{E}_{r}(x,A)$ is the value of the best approximation of
function  $x$ by trigonometric polynomials whose degree does not
exceed $r$ with respect to the norm in $\XX$. It is generally known
that differential operator $A$ is a generator of (isometric) group
of shifts in the space $\XX$:
\begin{eqnarray}
 & (U(t)x)(\xi)=x(t+\xi),\qquad x\in\XX;\,\, t,\xi\in\R\nonumber \\
 & \left\Vert U(t)\right\Vert \equiv1,\qquad t\in\R,\label{NormaZsyvy}\end{eqnarray}
where $\left\Vert U(\cdot)\right\Vert =\left\Vert
U(\cdot)\right\Vert _{\mathcal{L}(\XX)}$ is the norm of the operator
$U(t)$ in the space $\mathcal{L}(\XX)$ of linear continuous
operators over $\XX$. It follows from (\ref{NormaZsyvy}) that
\[
\widetilde{\omega}_{k}(t,x,A)=\omega_{k}(t,x,A)=\sup_{0\leq h\leq
t}\bigg\Vert \sum_{j=0}^{k}(-1)^{k-j}{{j \choose
k}}x(\cdot+jh)\bigg\Vert _{\XX},\quad t\in\R_{+},\, x\in\XX.
\]
I.e., in that case, $\widetilde{\omega}_{k}(t,x,A)$ coincides with
classic modulus of continuity of $k$-th degree in the space $\XX$.

Thus, from theorem \ref{Th1} and corollary \ref{NaslDlyaD(A^m)}
one can conclude all classic Jackson-type inequalities in the spaces
$C(2\pi)$ and $L_{p}(2\pi),$ $1\leq p<\infty$.
\end{example}

\subsection{Jackson's inequalities of the approximation by exponential type entire functions
 in the space $L_{p}(\R, \mu^p)$}
We consider the real-valued function $\mu(t)$ satisfying the following
conditions:
\begin{itemize}
\item[1)] $\mu(t)\ge 1,\quad t\in\R$;
\item[2)] $\mu(t)$ is even, monotonically non-decreasing when $t>0$;
\item[3)] $\mu(t)$ satisfies naturally occurring in many applications condition
 $\mu(t+s)\le \mu(t)\cdot\mu(s),\ s,t\in\R$.
\item[4)]
$\int_{-\infty}^{\infty}\frac{\ln\mu(t)}{1+t^2}\,dt<\infty$,
\end{itemize}
or alternatively, instead of 4), the equivalent condition holds:
\begin{itemize}
\item[4')] $\sum_{k=1}^\infty\frac{\ln\mu(k)}{k^2}<\infty$.
\end{itemize}

Lets consider several important classes of functions satisfying
conditions  1)--4).

1. Constant function $\mu(t)\equiv 1,\quad t\in\R$.

2. Functions with polynomial order of growth at infinity. It is easy
to check that for such functions following estimate holds: $\exists
k\in\N,\ \exists M\ge 1$
\begin{equation*}
\mu(t)\le M(1+|t|)^k,\quad t\in\R.
\end{equation*}

3. Functions of the form
\begin{equation*}
    \mu(t)=e^{|t|^\beta},\quad 0<\beta<1,\ t\in\R.
\end{equation*}

4. $\mu(t)$ represented as a power series for $t>0$. I.e.,
\begin{equation*}
    \mu(t)=\sum_{n=0}^\infty\frac{|t|^{n}}{m_n},
\end{equation*}
where $\{m_n\}_{n\in\N}$ is the sequence of positive real numbers
satisfying two conditions:
\begin{itemize}
\item $m_0=1$, $m_n^2\le m_{n-1}\cdot m_{n+1},\ n\in\mathbb{N}$;
\item $\forall k,l\in\mathbb{N}\ \frac{(k+l)!}{m_{k+l}}\le
\frac{k!}{m_k}\frac{l!}{m_l}$.
\end{itemize}
The function $\mu(t)$, defined above, obviously satisfies conditions
1) and 2). The condition $\forall k,l\in\mathbb{N}\
\frac{(k+l)!}{m_{k+l}}\le \frac{k!}{m_k}\frac{l!}{m_l}$ implies
\begin{equation}\label{example_inequal}
\sum_{k=0}^n\frac{t^ks^{n-k}n!}{k!(n-k)!\,m_n}\le\sum_{k=0}^n\frac{t^ks^{n-k}}{m_km_{n-k}},
\end{equation}
and it is easy to see that condition 3) follows from inequality
(\ref{example_inequal}). The Denjoy~- Carleman theorem
\cite[p.376]{Rudin1970} asserts that the following conditions are
equivalent:
\begin{itemize}
\item[a)] $\mu(t)$ satisfies condition 4);
\item[b)] $\sum_{n=1}^{\infty}\left(\frac
1{m_n}\right)^{1/n}<\infty$;
\item[c)] $\sum_{n=1}^{\infty}\frac{m_{n-1}}{m_n}<\infty$.
\end{itemize}

5. $\mu(t)$ as a module of an entire function with zeroes on the
imaginary axis. We consider
\begin{equation*}
    \omega(t)=C\prod_{k=1}^\infty\left(1-\frac
    t{it_k}\right),\quad t\in\R,
\end{equation*}
where $C\ge 1,\ 0<t_1\le t_2\le\ldots,\ \sum_{k=1}^\infty \frac
1{t_k}<\infty$. We set $\mu(t):=|\omega(t)|$. Then $\mu(t)$
satisfies conditions 1)~-- 3), and, as shown in \cite{LubMatsaev},
$\mu(t)$ satisfies condition 4) also.

Lets proceed to the description of the spaces $L_p(\R,\mu^p)$. Let
 the function $\mu(t)$ satisfies conditions 1)~-- 4). One can consider
 the space $L_p(\R,\mu^p)$ of the functions $x(s),\ s\in\R$, integrable in
$p$-th degree with the weight $\mu^p$:
\begin{equation*}
    \|x\|^p_{L_p(\R,\mu^p)}=\int_{-\infty}^{\infty}|x(s)|^p\mu^p(s)\,ds.
\end{equation*}
$L_p(\R,\mu^p)$ is the Banach space. We consider the differential
operator $A$ ($\D(A)=\{x\in L_p(\R,\mu^p)\cap AC(\R):\ x'\in
L_p(\R,\mu^p)\},\ (Ax)(t)=\frac {dx}{dt}$). As in example
\ref{Example2}, the operator $A$ generates the group of shifts
 $\{U(t)\}_{t\in\R}$ in the space $L_p(\R,\mu^p)$.
But in contrast to example \ref{Example2}, this group isn't
bounded. Indeed, lets consider
\begin{equation*}
    x(s)=\begin{cases}
    1, & s\in [0,1],\\
    0, & s\not\in [0,1].
    \end{cases}
\end{equation*}
Obviously, $x(s)\in L_p(\R,\mu^p)$, but for $t>1$
\begin{equation*}
    \|U(t)x\|^p=\int_{-\infty}^{\infty}|x(t+s)|^p\mu^p(s)\,ds=\int_{t-1}^t\mu^p(s)\,ds\ge\mu^p(t-1)\to\infty,\
    t\to\infty.
\end{equation*}

On the other hand, because of the property 3),
\begin{equation*}
    \|U(t)x\|^p=\int_{-\infty}^{\infty}|x(t+s)|^p\mu^p(s)\,ds\le
    \mu^p(-t)\int_{-\infty}^{\infty}|x(t+s)|^p\mu^p(t+s)\,ds=\big(\mu(-t)\big)^p\|x\|^p,
\end{equation*}
so $\|U(t)\|_{L_p(\R,\mu^p)}\le\mu(-t)=\mu(|t|),\quad t\in\R$.
 \footnote{If $\mu(t)$ is continuous and $\mu(0)=1$, it is possible to show in a similar manner that $\|U(t)\|_{L_p(\R,\mu^p)}=\mu(|t|)$.}

By the same way as in the example \ref{Example2}, modules of
continuity $\omega_k$ and $\tilde\omega_k$ coincides with classic
ones, but in contrast to the example \ref{Example2}, they don't
equal mutually. The space $\E(A)$ consists of fast decrescent at the
infinity entire functions. The examples of such functions have been
given in \cite{LubMatsaev}. By applying theorem \ref{Th1} one
can get
\begin{nasl}
$\forall k\in\N$ there exists constant $\mathbf{m}_{k}(p,\mu)>0$
such that $\forall f\in L_p(\R,\mu^p)$
\begin{equation*}
\mathcal{E}_{r}(f)\leq\mathbf{m}_{k}\cdot\tilde\omega_{k}\left(\frac{1}{r},x,A\right),\quad
r\ge 1.
\end{equation*}
\end{nasl}

\end{document}